\documentclass{article}

\usepackage{arxiv}

\usepackage[utf8]{inputenc} % allow utf-8 input
\usepackage[T1]{fontenc}    % use 8-bit T1 fonts
\usepackage{hyperref}       % hyperlinks
\usepackage{url}            % simple URL typesetting
\usepackage{booktabs}       % professional-quality tables
\usepackage{amsfonts}       % blackboard math symbols
\usepackage{nicefrac}       % compact symbols for 1/2, etc.
\usepackage{microtype}      % microtypography
\usepackage{lipsum}

\usepackage{graphicx} 
\usepackage{amsmath}
\usepackage{amssymb}
\usepackage{algorithm}
\usepackage{algorithmicx}
\usepackage{esint}% Include this line if your 

\newcommand{\C}{\mathbb{C}}

\newcommand{\Prob}[1]{\mathbb{P}\left( #1 \right)}

\if@secthm
  \newtheorem{thm}{Theorem}[section]
  \@addtoreset{thm}{section}
\else
  \newtheorem{thm}{Theorem}
\fi

\newtheorem{assum}[thm]{Assumption}

\newtheorem{defn}[thm]{Definition}

\newenvironment{pf}%
  {\par
  \noindent
   {\bfseries\Elproofname}\enspace\ignorespaces}%
  {\par
  }
\def\Elproofname{PROOF.}

\title{An overview of optimal control optimization problems driven by non-convexity measures}

\author{
  Weixin Wang
  Department of Industrial Engineering\\
  State University of New York\\
  Buffalo, NY 14260 \\
  \texttt{wangweixin23@gmail.com} \\
}

\begin{document}
\maketitle

\begin{abstract}
Recently, literature on dynamic coherent risk measures has broadened the choices for risk-sensitive performance evaluation. A running example includes Cumulative prospect theory and Conditional variance at risk. Most of them can be can be interpreted in general as a non-linear transformation of a given random variable. Non-convexity property has implied a lot of mathematical intricacies and challenges. The paper gives overview on the recent development of dynamic programming optimal control optimization problems driven by non-convex measures. 
\end{abstract}

% keywords can be removed
\keywords{Optimization \and Control \and Dynamic programming \and non-convexity}

\section{Introduction}
Since the introduction by Bellman \cite{Bellman1952}, Dynamic programming has
been the subject of extensive research in the past decades; see
for example \cite{Bellman1957}. Markov Decision Process(MDP) is a typical probabilistic model 
framework to study the Dynamic optimization problems.
In a MDP problem setting, outcomes are partly random and partly under the control of a decision maker.

In many applications, non-convex measures are more appropriate
than convex measures since risk-sensitive measures are better in capturing the real-world application case
\cite{Howard1972,Jaquette1973}.
In standard MDPs, the performance measures are frequently expressed
as expected utility functions that are risk-sensitive.
For example, many problems evaluate their outcomes by using $E\left[u(X)\right]$,
where $u$ is a risk-sensitive utility function (e.g., exponential),
and $X$ is a random variable representing the total reward or cost.
However, there are lot of mathematical intricacies implied by non-convex measures, and lots of theoretical properties of non-convex dynamic programming structure cannot be applied. 

The goal in a Markov decision process(MDP) or dynamic programming problem is to find a good "policy" for the decision maker. a function $\pi$  that specifies the action $\pi (s)$ that the decision maker will choose when in state $s$. Once a Markov decision process is combined with a policy in this way, this fixes the action for each state and the resulting combination behaves like a Markov chain (since the action chosen in state $s$ is completely determined by $\pi (s)$ and  $\Pr(s_{t+1}=s'\mid s_{t}=s,a_{t}=a)$ reduces to $\Pr(s_{t+1}=s'\mid s_{t}=s)$, a Markov transition matrix).

Convexity is a mathematical property that defining a real valued n-dimensional function such that line segment between any two points on the graph of the function lies above the graph between the two points. Equivalently, a function is convex if its epigraph (the set of points on or above the graph of the function) is a convex set. A twice-differentiable function of a single variable is convex if and only if its second derivative is nonnegative on its entire domain. Convex functions is important in many areas. They are especially important in the study of dynamic where they can derive a sequence of convenient properties. For example, contraction mapping or value iteration method can be justified to be converging into the maximal values. 

Many researchers have devoted their work on solving dynamic programming with non-convex measures. And this paper will give an overview on the models of dynamic programming under some special non-convex measures. One nominal example of non-convex measures is called Conditional value at risk(Cvar). Cvar is a concept used in the field of financial risk measurement to evaluate the market risk or credit risk of a portfolio. Cvar, is known to be a non-convex function and it is believed to lead a non-convex optimization problem in many settings. However, it is possible to transform the problem into a linear program and find the global solution. This property makes expected Cvar a cornerstone of alternatives to mean-variance portfolio optimization, which account for the higher moments (e.g., skewness and kurtosis) of a return distribution.

Another nominal example of non-convex measures is called cumulative prospect theory (CPT) 
developed by Tversky and Kahneman\cite{Tversky1992}. It is a more desired criteria to model human decision makers when the goal is to find the optimal policy whose outcome is maximally aligned with human's preference.
Unlike both expected utility and coherent risk
measures, which are normative approaches, CPT-based criteria have
risen from the search for a powerful descriptive model for human decision
making \cite{LIN20181}. Their ability to capture human decision dynamics under uncertainty
(e.g., lotteries) has strong empirical support \cite{Wakker2010}.
The incorporation of CPT into
dynamic systems has just been recently developed by \cite{LIN20181}. Meanwhile, He and Zhou \cite{He2011}
have studied a portfolio choice problem using a CPT-based approach.

Aside from CPT's ability to explain human decisions, the policies
it produces are randomized, which are more robust against modeling
errors\cite{LIN20181}. And that can be referred based on the 3 main properties of CPT:
First, it assumes that people tend to think of possible outcomes usually relative to a certain reference point (often the status quo) rather than to the final status, a phenomenon which is called framing effect. Moreover, people have different risk attitudes towards gains (i.e. outcomes above the reference point) and losses (i.e. outcomes below the reference point) and care generally more about potential losses than potential gains (loss aversion). Finally, people are usually a bit more optimistic about the chance of winning and tend to overweight extreme gains. 

A simple observation is that CPT is a generalization of expected utility so that convexity cannot be directly applied on it. Let $X$ be a Bernoulli random variable
that takes the value $1$ with probability $p$ and $0$ otherwise.
Since $E\left[u\left(X\right)\right]=pu(1)$, its expected utility
is always linear in $p$. Using a typical CPT weighting function $\frac{p^{\delta}}{\left(p^{\delta}+\left(1-p\right)^{\delta}\right)^{\frac{1}{\delta}}}$,
when $\delta=1$, the linear case is recovered. Convex risk measures
can also be recovered with appropriate weighting functions.

The paper is organized as follows. In section 2, we introduce cumulative
prospect theory and demonstrate the properties of CPT-based decisions.
In section 3, CPT-based criteria are applied to general dynamic problems.
In section 4, we are reviewing some recently development on the CPT-functional such as estimation and optimization.
In section 5, we are discussing some theorems recently developed on the CPT-functional driven dynamc programming problem. 

\section{Cumulative Prospect Theory: Mathematical formulation}

Prospect theory was built by Kahneman and Tversky
\cite{Kahneman1979}. The modified version, cumulative prospect theory(CPT) is developed in early 1990s
\cite{Tversky1992}. The 3 assumptions of CPT can be translated into the following mathematical expression:
1) The utility function has a reference point against which gains
and losses are evaluated, and this expression refers to the claim that people tend to think of possible outcomes usually relative to a certain reference point; 2) The utility function is concave on gains
and convex on losses, and this property just implies risk-avers behaviour of human being; 3) A probability
weighting function (cf. Def. \ref{def:probability-weighting-function})
that transforms the cumulative distribution function of a distribution such that
the probability of extreme is over weighted and the probability of common event is under weighted. We define the weighting function as the following:

\begin{defn}
\label{def:probability-weighting-function}A \textit{probability weighting
function}, $w$, is a monotonically non-decreasing continuous function
from $\left[0,1\right]$ to $[0,1]$ with $w\left(0\right)=0$ and
$w\left(1\right)=1$. 
\end{defn}

Let $X$ be a real random variable with a given probability distribution function, and $b$ to be denoted as a reference point together with two utilities functions $u^{+}$ and $u^{-}$. The $w^{+}$ and $w^{-}$ denote two different smooth probability weighting function. The CPT-functional applied on the random variable $X$ can be expressed as the following: 
\begin{flalign}
\mathbb{C}\left(X\right) & =\int_{0}^{\infty}w^{+}\left(\mathnormal{P}\left(u^{+}\left(\left(X-b\right)_{+}\right)>x\right)\right)dx\nonumber \\
 & -\int_{0}^{\infty}w^{-}\left(\mathnormal{P\left(u^{-}\left(\left(X-b\right)_{-}\right)>x\right)}\right)dx,\label{eq:CPT_main}
\end{flalign}
Notice that the notations $\left(\cdot\right)_{+}$
and $\left(\cdot\right)_{-}$ are shorthands for $\max\left(\cdot,0\right)$
and $\mbox{-\ensuremath{\min\left(\cdot,0\right)}}$, respectively.
Appropriate integrability assumptions are satisfied. 

\section{Dynamic Programming driven by CPT-functional}
\subsection{Simple example on expected utility}
This part of the chapter follows closely the instruction of abstract dynamic programming from the chapter of 
\cite{Bertsekas2013}. In an abstract dynamic programming framework, the key components include a state space $\mathbb{X}$, an action space $\mathbb{A}$ and a disturbing noise space $\delta_{0}$ on the real valued space $\mathbb{R}$. 

The system of the dynamic programming evolves according to a specific dynamic $x_{t+1}=f(x_{t},a_{t},\delta_{t}),$. In the above formula, $x_{t+1}$ and $x_{t}$ belong to the state space $\mathbb{X}$, and action belongs to the action space $\mathbb{A}$ and the disturbing noise $\delta_{t}$ has a specific distribution $\mathcal{P}_t$. 

At each step $t$, the combination of action, state and noise will together determine a reward $g(x_{t},a_{t}, \delta_{t})$ with g defined as a mapping from the space $\mathbb{X} \times \mathbb{A} \times \mathbb{R}$ towards real valued space $\mathbb{R}$. 

In the may classical cases, at each time stamp, we are interested in the expected return of the reward
$E\left[g(x_{0},a_{0},\delta_{0}) \right]$. Since the value at the previous states may have a impact on the return in the future states. 
The standard expected value across the time steps can be written as
$$ E\left[g(x_{0},a_{0},\delta_{0})+
E\left[g\left(x_{1},a_{1},\delta_{1}\right)+\cdots|x_{1}\right]|x_{0}\right]$$
.

\subsection{Dynamic programming: Pure abstract structure}

We are interested in nonempty Borel spaces $X$ and $A$ of states
and controls such that for each $x\in X$ there is a nonempty feasible
control Borel set $A\left(x\right)\subset A$. We denote the set of
probability measures over $A$ equipped with the Prohorov metric by
$\mathcal{P}\left(A\right)$. We denote by $\mathcal{S}$ the set
of all measurable functions $\mu:X\rightarrow\mathcal{P}\left(A\right)$
satisfying $\mu\left(x\right)\in\mathcal{P}\left(A\left(x\right)\right),\;\forall x\in X$,
which we refer to as policies. The nonempty Borel space of disturbances
is denoted by $\Delta,$ and given a state-action pair $\left(x_{k},a_{k}\right)\in X\times A$,
an element $\delta_{k}\in\Delta\left(x_{k},a_{k}\right)\subset\Delta$
drives the system to its next state through a measurable function
$f:X\times A\times\Delta\rightarrow X$ by $x_{k+1}=f\left(x_{k},a_{k},\delta_{k}\right)$.
At each time $k$, a per-step cost is accumulated and denoted by a
measurable function $g:X\times A\times\Delta\rightarrow\mathbf{R}$.
The stochastic kernel $\mathnormal{P}\left(\cdot|x,a\right)$ is defined
over $\Delta\left(x,a\right)$. Furthermore, we denote both the realization
and the random variable disturbance at time $k$ by $\delta_{k}$.
We denote by $R\left(X\right)$ the set of real-valued measurable
functions $J:X\rightarrow\mathbf{R}$. A nonstationary Markov policy
is denoted by $\pi=\left\{ \mu_{0},\mu_{1},\mu_{2},\dots\right\} ,$
where $\mu_{k}\in\mathcal{S}$ and $\Pi$ denotes the set of all feasible
non-stationary Markov policies. When the context is clear, we will
refer to $\pi$ as policies as well. 

Given an element $\bar{J}\in R\left(X\right)$, we minimize the cost
over all non-stationary Markov policies, i.e., 
\begin{eqnarray}
J^{*}\left(x\right) & = & \inf_{\pi\in\Pi}J_{\pi}\left(x\right),\text{ where}\nonumber \\
J_{\pi}\left(x\right) & = & \limsup_{k\rightarrow\infty}\left(T_{\mu_{0}}T_{\mu_{1}}T_{\mu_{2}}\cdots T_{\mu_{k}}\bar{J}\right)\left(x\right),\label{eq:main-optimization-problem}
\end{eqnarray}

for all $x\in X,$ and $T_{\mu}:R\left(X\right)\rightarrow R\left(X\right)$
is a problem dependent operator. 
We define a mapping: 
$$H:X\times\mathcal{P}\left(A\right)\times R\left(X\right)\rightarrow\mathbf{R}$$
such that for each policy $\mu\in\mathcal{S}$ it satisfies $\left(T_{\mu}J\right)\left(x\right)=H\left(x,\mu\left(x\right),J\right),\;\forall x\in X.$
We define the operator $T$ by 
$$\left(TJ\right)\left(x\right)=\inf_{a\in\mathcal{P}\left(A\left(x\right)\right)}H\left(x,a,J\right)=\inf_{\mu\in\mathcal{S}}\left(T_{\mu}J\right)\left(x\right),\;\forall x\in X.$$ 

There exist two infinite horizon problems, namely discounted and
transient, using Eq. \ref{eq:main-optimization-problem}. The approach
is to apply Bertsekas's abstract dynamic programming \cite{Bertsekas2013}
to both problems. In both cases, we are looking for assumptions to
satisfy the \emph{Monotonicity} and contraction assumptions in his approach.
Those assumptions can lead to the conclusion that value
and policy iteration converge to a unique value function and an
optimal policy can be attained (\cite{Bertsekas2013}). Additionally such value function attained by the optimal policy can be approximated within arbitrary accuracy by a stationary policy.  
The two assumptions are \emph{Monotonicity} and \emph{Contraction}. 
(see \cite{Bertsekas2013}, Assumptions 2.1.1 and 2.1.2).

\begin{assum}\label{assum:monotonicty}

(Monotonicity) If $J,J'\in R\left(X\right)$ and $J\leq J'$, then
$H\left(x,a,J\right)\leq H(x,a,J'),\;\forall x\in X,\ a\in\mathcal{P}\left(A\left(x\right)\right).$ 

\end{assum}

\begin{assum}\label{assum:contraction}
(Contraction) For all $J\in B\left(X\right)$ and $\mu\in\mathcal{S}$,
the functions $T_{\mu}J$ and $TJ$ belong to $B\left(X\right)$.
Furthermore, for some $\alpha\in\left(0,1\right),$ we have $\left\Vert T_{\mu}J-T_{\mu}J'\right\Vert \leq\alpha\left\Vert J-J'\right\Vert ,\;\forall J,J'\in B\left(X\right),\ \mu\in\mathcal{S}$. 
\end{assum}

Moreover, we assume that function $g$ satisfies the following
assumption: 

\begin{assum}\label{assum:per-stage-boundedness}

There exists a constant $c>0$ such that $\sup_{x\in X,a\in A\left(x\right),\delta\in\Delta}\left|g\left(x,a,\delta\right)\right|\leq c.$ 

\end{assum}

Assumptions \ref{assum:monotonicty} and \ref{assum:contraction}
are conditions on the $H$. And assumption \ref{assum:per-stage-boundedness}
is standard for the classical discounted infinite horizon problem. 

\subsection{Definition of $H$ under CPT functional}

Assume that we are given a system with dynamics $x_{k+1}=f\left(x_{k},a_{k},\delta_{k}\right),$
with a discount factor of $\alpha\in\left(0,1\right)$. If we replace the expected utility function above by the CPT-functional, the corresponding
$H$ mapping for the problem is $H\left(x,a,J\right)=$
\begin{gather}
\int_{0}^{\infty}w_{+}\left(\mathnormal{P}\left(u_{+}\left(\left(g\left(x,a,\delta\right)+\alpha J\left(f\left(x,a,\delta\right)\right)\right)_{+}\right)>z\right)\right)dz\nonumber \\
\mbox{-}\int_{0}^{\infty}w_{-}\left(\mathnormal{P}\left(u_{-}\left(\left(g\left(x,a,\delta\right)+\alpha J\left(f\left(x,a,\delta\right)\right)\right)_{-}\right)>z\right)\right)dz,\label{eq:Discounted-Infinite-Horizon}
\end{gather}

Without loss
of generality, the reference point 
is assumed to be zero for Eq. \ref{eq:Discounted-Infinite-Horizon}. 

The work from \cite{LIN20181} has proved the following theorem:
\begin{thm}
\label{thm:monotonicity-discounted}If $u_{+},u_{-}:\mathbf{R}^{+}\rightarrow\mathbf{R}^{+}$
are both monotonically non-decreasing functions, and $w_{+},w_{-}$
are probability weighting functions, then Eq. \ref{eq:Discounted-Infinite-Horizon}
satisfies Assumption \ref{assum:monotonicty}. \end{thm}

\begin{pf}
Refer to \cite{LIN20181}. 
\end{pf}

Next, the following theorem proved by \ref{eq:Discounted-Infinite-Horizon} states that $T_{\mu}$ defined by $H$ is a contraction. 
\begin{thm}
\label{thm:contraction-discounted}Assume the following conditions
hold: 1) the assumptions in Theorem \ref{thm:monotonicity-discounted}
hold; 2) $u_{+},u_{-}$ are invertible (denoted by $u_{+}^{-1}$ and
$u_{-}^{-1}$), differentiable (denoted by $u'_{+}$ and $u'_{-}$)
with $u_{+}\left(0\right)=u_{-}\left(0\right)=0$; 3) $u'_{+},u'_{-}$
are monotonically non-increasing; 4) there exists a $\beta\in\left(0,1\right)$
such that the inequality $\int_{0}^{\alpha c}w_{+}\left(\mathnormal{P}\left(Z<z\right)\right)u'_{+}\left(\alpha c-z\right)dz+\int_{0}^{\alpha c}w_{-}\left(\mathnormal{P}\left(Z>z\right)\right)u'_{-}\left(z\right)dz\leq\beta c,\ c>0$
holds for any non-negative real-valued random variable $Z$. Then
Assumption \ref{assum:contraction} is satisfied. \end{thm}
%\vspace{-3cm}

\begin{pf}
Refer to \cite{LIN20181}
\end{pf}

\subsection{Transient Markov Control Model }

%In this section, we prove the optimality of the dynamic programming
% equation for transient Markov control models. 
A transient Markov model
evolves according the equation $x_{t+1}=f\left(x_{t},a_{t},\delta_{t}\right)$
and has some absorbing state $x_{A}\in X$, such that if $x_{t}=x_{A}$,
then $f(x_{A},a,\delta)=x_{A}$ and $g\left(x_{A},a,\delta\right)=0$
for all $a\in A\left(x_{A}\right),\ \delta\in\Delta$. To interpret the formula, note that once we reach absorbing state, the episode of a dynamic programming is terminated and no more action needs to be taken. We denote the first
hitting time of the absorbing state with a policy $\pi\in\Pi$ by
$\tau_{A}^{\pi}:=\inf\left\{ t\geq0|x_{t}^{\pi}=x_{A}\right\} $.
Typically we can assume that a transient Markov model reaches its absorbing state in a finite amount
of time starting from an initial state $x_{0}$.

In the absorbing transient Markov Control model, the corresponding $H$ mapping for the systems is: $H\left(x,a,J\right)=$
\begin{gather}
\int_{0}^{\infty}w_{+}\left(\mathnormal{\tilde{P}}\left(u_{+}\left(\left(g\left(x,a,\delta\right)+J\left(f\left(x,a,\delta\right)\right)\right)_{+}\right)>z\right)\right)dz\nonumber \\
-\int_{0}^{\infty}w_{-}\left(\mathnormal{\tilde{P}}\left(u_{-}\left(\left(g\left(x,a,\delta\right)+J\left(f\left(x,a,\delta\right)\right)\right)_{-}\right)>z\right)\right)dz.\label{eq:tranient-cannonical}
\end{gather}
where $\tilde{P}$ is defined as $\tilde{P}\left(\cdot\right)=P\left(\cdot\cap f\left(x,a,\delta\right)\in X\setminus x_{A}\right)\leq1,\;\forall a\in A\left(x\right),\ x\in X$.
\begin{defn}
A policy $\pi=\left\{ \mu_{0},\mu_{1},\dots\right\} \in\Pi$ is \textit{transient}
with respect to a Markov control model, if there exists a constant
$c$ such that $\sum_{k=0}^{\infty}\mathnormal{P}\left(f\left(x_{k},\mu_{k}\left(x_{k}\right),\delta_{k}\right)\in X\setminus x_{A}\right)\leq c.$
If the inequality above holds for all $\pi\in\Pi$, then the model
is called \textit{uniformly transient. }The inequality is also known
as the Pliska condition \cite{Pliska1978}.
\end{defn}

The next theorem gives the conditions needed to satisfy the monotonicity
and contraction assumptions. 
\begin{thm}
If $u_{+},u_{-}:\mathbf{R}^{+}\rightarrow\mathbf{R}^{+}$ are both
monotonically non-decreasing functions, and $w_{+},w_{-}$ are probability
weighting functions, then Eq. \ref{eq:tranient-cannonical} satisfies
Assumption \ref{assum:monotonicty}.\end{thm}

\begin{pf}
Use the same argument as in Theorem \ref{thm:monotonicity-discounted}.\end{pf}

\begin{thm}
\label{thm:Transient-Markov-Model-Three-Conditions-Are-Satisfied}Assume
the following conditions hold:  1) the Markov control model is uniformly
transient; 2) $\exists C>0$ such that $u'_{+}\left(0\right)\leq C,\ u'_{-}\left(0\right)\leq C$;
3) conditions 1-3 in Theorem \ref{thm:contraction-discounted} hold;
4) $\exists\xi>0$ such that $w_{+}\left(x\right)\leq\xi x$ and $w_{-}\left(x\right)\leq\xi x.$
Then the operator $T_{\mu}$ defined by using Eq. \ref{eq:tranient-cannonical}
is a $K$-step contraction. \end{thm}

\begin{pf}
Refer to \cite{LIN20181}
\end{pf}

\section{Numerical development on CPT-functional}
This chapter discusses the recent development in the numerical study of CPT-functional. The paper  \cite{la2016cumulative} developed a quantile-based statistical estimator on the CPT-functional applied on a given random variable, and the paper \cite{8329994} also developed a few stochastic optimization techniques to optimize the policy in a dynamic programming setting given CPT-return.

\subsection{CPT-functional estimator}
We are referring to the quantile-based estimator developed by \cite{la2016cumulative}. Recall the CPT-functional defined by the formula \ref{eq:CPT_main}, to make it simple and without loss of generality, we assume the benchmark $b = 0$ and we are only focusing on estimating the positive part of CPT-functional $\int_0^{+\infty} w^+\left(\Prob{u^+(X)>x}\right) dx$. The following estimation scheme \ref{alg:holder-est} developed by \cite{la2016cumulative} is providing an asymptotic consistent estimator of the CPT-functional.

\begin{defn} \textbf{CPT-value estimation}
\label{alg:holder-est}
Assume that we simulate $n$ i.i.d. samples from the distribution of $X$. We denote the order statistics of the the samples and label them as follows: 
$X_{[1]}, X_{[2]}, \ldots ,X_{[n]}$. Note that $u^+(X_{[1]}),\ldots ,u^+(X_{[n]})$ are also in ascending order.
Let
\vspace{-0.5ex}
$$\hspace{-1.1em}\overline \C_n^+:=\sum_{i=1}^{n} u^+(X_{[i]}) \left(w^+\!\left(\frac{n+1-i}{n}\right)\!-\! w^+\!\left(\frac{n-i}{n}\right) \right).$$
\vspace{-0.5ex}

Apply $u^{-}$ on the sequence $\{X_{[1]}, X_{[2]}, \ldots ,X_{[n]}\}$; notice that $u^{-}(X_{[i]})$ is in descending order since $u^{-}$ is a decreasing function.     

Let
\vspace{-0.5ex}
$$\overline \C_n^-:=\sum_{i=1}^{n} u^-(X_{[i]}) \left(w^-\left(\frac{i}{n}\right)- w^-\left(\frac{i-1}{n}\right) \right). $$
\vspace{-0.5ex}
The  $\overline \C_n =\overline \C_n^+ - \overline \C_n^-$ is an estimator of the CPT-functional defined in the formula \ref{eq:CPT_main}
\end{defn}

The theorem below proved the asymptotic consistency of CPT-functional

\begin{thm}(\textbf{Asymptotic consistency})
\label{prop:holder-asymptotic}
Assume (A1) and that $F^+(\cdot)$ and $F^-(\cdot)$, the respective distribution functions of $u^+(X)$ and $u^-(X)$, 
are Lipschitz continuous 
%with constants $L^+$ and $L^-$, respectively, 
on the respective intervals $(0,+\infty)$, and 
$(-\infty, 0)$. Then, we have that
\begin{align}
\overline \C_n
\rightarrow
\C(X)
 \text{   a.s. as } n\rightarrow \infty
\end{align}
where $\overline \C_n$ is as defined in Algorithm \ref{alg:holder-est} and $\C(X)$ as in \eqref{eq:CPT_main}.
\end{thm}

\begin{pf}
Refer to \cite{la2016cumulative}
\end{pf}

\bibliographystyle{unsrt}  
\bibliography{references} 

\end{document}